\documentclass[12pt, reqno]{amsart}
\usepackage{amsmath}
\usepackage{amssymb, amsfonts}

\numberwithin{equation}{section}

\begin{document}

{\theoremstyle{theorem}
    \newtheorem{theorem}{\bf Theorem}[section]
    \newtheorem{proposition}[theorem]{\bf Proposition}
    \newtheorem{claim}[theorem]{\bf Claim}
    \newtheorem{lemma}[theorem]{\bf Lemma}
    \newtheorem{corollary}[theorem]{\bf Corollary}
}
{\theoremstyle{remark}
    \newtheorem{remark}[theorem]{\bf Remark}
    \newtheorem{example}[theorem]{\bf Example}
}
{\theoremstyle{definition}
    \newtheorem{defn}[theorem]{\bf Definition}
    \newtheorem{question}[theorem]{\bf Question}
}

%%%%%%%%%%%%%%%%%%%%%%%%%%%%%%%%%%%%%%%%%%%

\newcommand{\lreg}[1]{{\operatorname{reg}_{\mathfrak L}^{[#1]}}}
\newcommand{\la}[1]{a^{[#1]}}
\newcommand{\rreg}[1]{{{\operatorname{reg}^{\mathfrak V}_{#1}}}}
\newcommand{\ppi}[1]{\pi^{[#1]}}
\newcommand{\X}[1]{{X^{[#1]}}}
\newcommand{\mS}[1]{{S^{[#1]}}}
\newcommand{\lF}[1]{F^{[#1]}}
\newcommand{\kset}[1]{[#1]}
\renewcommand{\c}[1]{c^{[#1]}}
\newcommand{\lr}[1]{r_{#1}}
\newcommand{\supp}[1]{{\operatorname{Supp}^{[#1]}}}
\newcommand{\s}{\bigskip}
\newcommand{\ms}{\medskip}
\newcommand{\sms}{\smallskip}
\newcommand{\f}[1]{\ensuremath{\mathfrak{#1}}}
\newcommand{\rom}[1]{\ensuremath{\textrm{#1}}}

\def\Ass{\operatorname{Ass}}
\def\Ann{\operatorname{Ann}}
\def\codim{\operatorname{codim}}
\def\coker{\operatorname{coker}}
\def\depth{\operatorname{depth}}
\def\deg{\operatorname{deg}}
\def\dim{\operatorname{dim}}
\def\Ext{\operatorname{Ext}}
\def\gr{\operatorname{gr}}
\def\grade{\operatorname{grade}}
\def\Hom{\operatorname{Hom}}
\def\height{\operatorname{ht}}
\def\id{\operatorname{id}}
\def\im{\operatorname{im}}
\def\ker{\operatorname{ker}}
\def\La{{\mathfrak a}_{\mathfrak L}}
\def\Lreg{\operatorname{\bf reg}^{\mathfrak L}}
\def\Max{\operatorname{Max}}
\def\Min{\operatorname{Min}}
\def\Proj{\operatorname{Proj}}
\def\projdim{\operatorname{projdim}}
\def\Quot{\operatorname{Quot}}
\def\rank{\operatorname{rank}}
\def\reg{\operatorname{reg}}
\def\Reg{\operatorname{\bf reg}}
\def\Rreg{\operatorname{\bf reg}^{\mathfrak V}}
\def\sign{\operatorname{sgn}}
\def\Spec{\operatorname{Spec}}
\def\Supp{\operatorname{Supp}}
\def\to{\longrightarrow}
\def\To{{\longrightarrow}}
\def\Tor{\operatorname{Tor}}
\def\V{\operatorname{V}}

\def\A{{\mathcal A}}
\def\D{\mathcal{D}}
\def\F{{\mathcal F}}
\def\G{{\mathcal G}}
\def\H{{\mathcal H}}
\def\I{{\mathcal I}}
\def\L{{\mathcal L}}
\def\M{{\mathcal M}}
\def\N{{\mathcal N}}
\def\O{{\mathcal O}}
\def\R{{\mathcal R}}
\def\T{{\mathcal T}}
\def\U{{\mathfrak U}}
\def\GG{{\mathbb G}}
\def\NN{{\mathbb N}}
\def\PP{{\mathbb P}}
\def\ZZ{{\mathbb Z}}
\def\VV{{\mathbb V}}
\def\a{{\mathbf a}}
\def\mba{{\bf a}}
\def\d{{\bf d}}
\def\e{{\bf e}}
\def\f{{\bf f}}
\def\k{{\bf k}}
\def\m{{\bf m}}
\def\n{{\bf n}}
\def\p{{\bf p}}
\def\q{{\bf q}}
\def\x{{\bf x}}
\def\y{{\bf y}}
\def\v{{\bf v}}
\def\1{{\bf 1}}
\def\0{{\bf 0}}
\def\aa{{\mathfrak a}}
\def\gm{{\mathfrak m}}
\def\pp{{\mathfrak p}}
\def\qq{{\mathfrak q}}
\def\r{{\mathbf r}^{\mathfrak L}}
\def\dd{{\mathfrak d}}
\def\fM{{\mathfrak M}}
\def\smap{{\longrightarrow\!\!\!\!\!\!\!\!\!\longrightarrow}}

\def\CM{Cohen-Macaulay}
\def\CMs{Cohen-Macaulay }

\long\def\alert#1{\smallskip{\hskip\parindent\vrule%
\vbox{\advance\hsize-2\parindent\hrule\smallskip\parindent.4\parindent%
\narrower\noindent#1\smallskip\hrule}\vrule}}

\title{Cohen-Macaulay multigraded modules}

\author{C-Y. Jean Chan}
\address{Department of Mathematics, University of Arkansas, Fayetteville AR 72701, USA}
\email{cchan@uark.edu}
\urladdr{http://comp.uark.edu/~cchan}

\author{Christine Cumming}
\address{Department of Mathematics and Physics, University of Louisiana at Monroe, Monroe, LA 71209, USA}
\email{cumming@ulm.edu}
\urladdr{http://www.ulm.edu/~cumming}

\author{Huy T\`ai H\`a}
\address{Tulane University, Department of Mathematics, 6823 St. Charles Ave., New Orleans LA 70118, USA}
\email{tha@tulane.edu}
\urladdr{http://www.math.tulane.edu/~tai/}

\subjclass[2000]{13A30, 13C14, 14B15, 14M05}
\keywords{Cohen-Macaulay, multigraded module, Rees algebra, Rees module, multi-Rees algebra, multi-Rees module}

\begin{abstract} Let $S$ be a standard $\NN^r$-graded algebra over a local ring $A$, and let $M$ be a finitely generated $\ZZ^r$-graded module over $S$. We characterize the Cohen-Macaulayness of $M$ in terms of the vanishing of certain sheaf cohomology modules. As a consequence, we apply our result to study the Cohen-Macaulayness of multi-Rees modules. Our work extends previous studies on the Cohen-Macaulayness of multi-Rees algebras.
\end{abstract}

\maketitle

%%%%%%%%%%%%%%%%%%%%%%%%%%%%%%%%%%%%%%%%%%%%%%%%%%%%%%%%%%%%%%%%%%%%%%

\section{Introduction} \label{s.intro}

The notion of Cohen-Macaulay rings and modules marks the interplay
between powerful lines of research in commutative algebra, algebraic
geometry, and algebraic combinatorics. It finds surprising
applications in far reaching problems and topics, for instance, in
duality theory, in homological theory of rings, and in the study of
polytopes and simplicial complexes.

Let $(A,\gm)$ be a local ring. Let $I \subseteq A$ be a proper
ideal, and let $\R = A \oplus It \oplus I^2t^2 \oplus \dots \subset
A[t]$ be the Rees algebra of $I$. Besides encoding many algebraic
properties of the ideal $I$ as well as its powers, the Rees algebra
$\R$ also gives an algebraic realization of the blowing up of $\Spec
A$ at the subscheme defined by $I$. Thus, characterizing the
Cohen-Macaulayness of $\R$ has always been an important problem in
commutative algebra. Lipman \cite{LI} succeeded in using Sancho de
Salas sequences to study the Cohen-Macaulayness of $\R$ via the
vanishing of sheaf cohomology groups on the blowup $\Proj \R$.

In recent years, much effort has been put forward to extend our
knowledge from the $\ZZ$-graded case to a more general multi-graded
setting (cf. \cite{HHR1,HHR2,HY,JV1,JV2,KN,V}). Lipman's method was generalized by Hyry \cite{HY} to
investigate the Cohen-Macaulayness of standard $\NN^r$-graded
algebras over a local ring. More precisely, \cite[Theorem 3.1]{HY}
shows that if $S = \bigoplus_{\n \ge \0}S_\n$ is a standard
$\NN^r$-graded algebra over $(A,\gm)$ such that its irrelevant ideal
$S_+ = \bigoplus_{\n > \0} S_\n$ has positive height, $Z = \Proj S$,
and $E = Z \times_A A/\gm$, then $S$ is a Cohen-Macaulay ring with a
negative $a$-invariant $\mba(S) < \0$ if and only if the following
conditions are satisfied:
\begin{itemize}
\item $\Gamma(Z, \O_Z(\n)) = S_\n$ for all $\n \ge \0$,
\item $H^i(Z, \O_Z(\n)) = 0$ for all $i > 0$ and $\n \ge \0$,
\item $H^i_E(Z, \O_Z(\n)) = 0$ for all $i < \dim Z$ and $\n < \0$.
\end{itemize}

The goal of this paper is to extend Hyry's result to study the
Cohen-Macaulayness of arbitrary finitely generated $\ZZ^r$-graded
modules over $S$. Let $M$ be a finitely generated $\ZZ^r$-graded
$S$-module, and let $\M$ be its associated coherent sheaf on $Z$.
Our first result, Theorem \ref{mainthm}, gives a
characterization for the Cohen-Macaulayness of $M$ in terms of the
vanishing of sheaf cohomology groups of twisted modules $\M(\n)$ on
$Z$ and with support $E$.

We also apply Theorem \ref{mainthm} to study the
Cohen-Macaulayness of multi-Rees modules. It is well-known (cf.
\cite{HHR1, HHR2, HY}) that if $I_1, \dots, I_r \subset A$ are
ideals of positive heights such that the multi-Rees algebra of $I_1,
\dots, I_r$ is Cohen-Macaulay then the usual Rees algebra of the
product $I_1 \cdots I_r$ is also Cohen-Macaulay. Our next
result, Theorem \ref{Rees module}, extends this phenomenon to
multi-Rees modules.

The converse of Theorem \ref{Rees module}, even in the case of
multi-Rees algebras, is known to be false. It is then desirable to
seek for conditions which, together with the Cohen-Macaulayness of
the Rees module of $I_1 \cdots I_r$ with respect to a given
$A$-module $N$, would imply that the multi-Rees module of $I_1,
\dots, I_r$ with respect to $N$ is Cohen-Macaulay. Hyry \cite{HY}
solved this problem for multi-Rees algebras (i.e. when $N = A$)
provided that the analytic spread of $I_1 \cdots I_r$ is small. Our
last result, Theorem \ref{Rees converse}, shows that the general
problem for multi-Rees modules, under some additional conditions, has a similar solution.

To prove Theorem \ref{mainthm}, we investigate local cohomology of
the Rees module of the irrelevant ideal $S_+$ with respect to $M$
under various graded structures. Here is a summary of the main ideas
of the proof. Let $R = \R_S(S_+)$ and $T = \R_M(S_+)$ be the Rees
algebra and Rees module of $S_+$ with respect to $M$, respectively.
Clearly, $R$ is an $\NN$-graded algebra over $S$ and $T$ is a
finitely generated $\ZZ$-graded $R$-module. The ring $S$ can also be
viewed as a standard $\NN$-graded algebra over $A$ (by coarsening
the graded structure). Let $\fM_R$ and $\fM_S$ be the maximal
homogeneous ideals in $R$ and in $S$, respectively. Let $Y = \Proj
R$ and $F = Y \times_S S/\fM_S$.  Let $\T$ be the associated
coherent sheaf of $T$ on $Y$. At the heart of our arguments is the
following Sancho de Salas sequence (cf. \cite[p. 150]{LI})
\begin{align}
\cdots \rightarrow [H^i_{\fM_R}(T)]_0 \rightarrow H^i_{\fM_S}(M)
\rightarrow H^i_F(Y, \T) \rightarrow [H^{i+1}_{\fM_R}(T)]_0
\rightarrow \cdots \label{SSsequence}
\end{align}

We start by observing that $R$ and $T$ have a natural $\ZZ^{r+1}$-graded structure
given by
$$R = \bigoplus_{\n \in \NN^r, k \ge 0} R_{(\n;k)} \mbox{ and
} T = \bigoplus_{\n \in \ZZ^r, k \ge 0} T_{(\n;k)},$$ where
$R_{(\n;k)} = S_{(n_1+k, \dots, n_r+k)}t^k$ and $T_{(\n;k)} = M_\n
[S_+^k]_{(k, \dots, k)}t^k$ for $\n = (n_1, \dots, n_r)$. The
cohomology modules $H^i_{\fM_R}(T)$ then inherit this
$\ZZ^{r+1}$-graded structure, and we can write $$[H^i_{\fM_R}(T)]_0
= \bigoplus_{\n \in \ZZ^r} H^i_{\fM_R}(T)_{(\n;0)}.$$ Next, we show
that $H^i_{\fM_R}(T)_{(\n;k)} = 0$ for $k \ge 0$ and $\n < \v(M)$.
This is done in Lemma \ref{L:3.3}. Together with the sequence
(\ref{SSsequence}), this implies that $[H^i_{\fM_S}(M)]_\n =
[H^i_F(Y,\T)]_\n$ for $\n < \v(M)$. We now observe that the
Cohen-Macaulayness of $M$ is characterized by the vanishing of
$H^i_{\fM_S}(M)$ for $i < \dim M$. Theorem \ref{mainthm} is then
proved by establishing the relationship between $H^i_F(Y,\T)$ and
$H^{i-r}_E(Z,\M)$ and the vanishing of $[H^i_{\fM_S}(M)]_\n$ for $\n
\not< \v(M)$. These are done in Lemma \ref{P:3.1}.

We start our proof of Theorem \ref{Rees module} by showing that if
$M$ is the multi-Rees module of $I_1, \dots, I_r$ with respect to an
$A$-module $N$ then the $a$-invariant of $M$ can be calculated
explicitly, namely $\mba(M) = - \1$. This is done in Lemma
\ref{L:2.1}. Observe further that $\v(M) = \0 > -\1$ in this case,
and so Theorem \ref{mainthm} can be applied. Next, we let $S$ be the
multi-Rees algebra of $I_1, \dots, I_r$, then the Rees algebra of
the product $I_1 \cdots I_r$ is a diagonal subalgebra $S^\Delta$ of
$S$ (which is $\NN$-graded). Theorem \ref{Rees module} is now proved
by noticing that there is a canonical isomorphism $f: \Proj S
\longrightarrow \Proj S^\Delta$ and pushing forward through $f$ to
reduce the problem to the well known $\ZZ$-graded situation.

Our last theorem, Theorem \ref{Rees converse}, is proved by a
straightforward generalization of Hyry's method in \cite{HY} from
multi-Rees algebras to multi-Rees modules. The paper is outlined as
follows. In Section \ref{s.prel}, we collect the notation, the
terminology, and the basic results that will be used throughout the
paper. Section~\ref{s.main} is devoted to proving the main theorem,
Theorem~\ref{mainthm}, that characterizes the Cohen-Macaulayness of
a finitely generated multi-graded module using sheaf cohomology
modules. Finally in Section~\ref{s.applications}, as an application,
we further deduce conditions to when the Cohen-Macaulayness of
multi-Rees modules of ideals $I_1, \cdots, I_r$ with respect to a
module and that of its diagonal submodule become equivalent.

\noindent{\bf Acknowledgment:} This project started when Hurricane
Katrina was about to hit New Orleans, where the last two authors
were at the time. The collaboration was then made possible with the
help of many people. In particular, we would like to thank Ian M.
Aberbach, S. Dale Cutkosky, Mark Johnson, Hema Srinivasan and Bernd
Ulrich. We would also like to thank the University of
Missouri-Columbia and the University of Arkansas for their
hospitality. The last author is partially supported by Louisiana Board of Regents Enhancement Grant.

%%%%%%%%%%%%%%%%%%%%%%%%%%%%%%%%%%%%%%%%%%%%%%%%%%%%%%%%%%%%%%%%%%%%%%

\section{Preliminaries} \label{s.prel}

% The theory of multi-projective schemes, multigraded rings and
% modules used in this paper is similar to that of projective schemes,
% graded rings and modules. In fact, any multi-projective scheme is
% projective by the isomorphism induced by a sufficiently high degree
% diagonal embedding.

For elementary facts about schemes, graded rings, and local cohomology
modules, we refer the reader to \cite{bs, EGA2, EGA3, har}.

Let $\0 = (0, \dots, 0)$ and $\1 = (1, \dots, 1)$. Let $\e_1, \dots,
\e_r$ be the standard basis vectors of $\ZZ^r$. Throughout the
paper, $S = \bigoplus_{\n \ge \0} S_\n$ will denote a standard
$\NN^r$-graded algebra over a local ring $(A,\gm)$. That is, $S$ is
generated over $S_\0 = A$ by elements of $\bigoplus_{j=1}^r
S_{\e_j}$. Define $S_+$ to be the {\it irrelevant ideal} of $S$
which is $\bigoplus_{\n > \0} S_\n$.  Let $S^\Delta = \bigoplus_{n
\ge 0} S_{(n, \dots, n)}$ denote the {\it diagonal subalgebra} of
$S$. Also, $M = \bigoplus_{\n \in \ZZ^r} M_\n$ will denote a
finitely generated $\ZZ^r$-graded $S$-module. Set $M^\Delta =
\bigoplus_{n \in \ZZ} M_{(n, \dots, n)}$.  We will call $M^\Delta$
the {\it diagonal submodule} of $M$. Clearly, $M^\Delta$ is a
$\ZZ$-graded $S^\Delta$-module.

For a vector $\n \in \ZZ^r$, we always use $n_1, \dots, n_r$ to represent its coordinates. For $\n, \m \in \ZZ^r$, we shall write $\n \ge \m$ if $n_j \ge m_j$ for all $j = 1, \dots, r$;
similarly, we write $\n > \m$ if $n_j > m_j$ for all $j = 1, \dots,
r$. We also define
\begin{align*}
& \min \{ \n, \m \} = (\min\{n_1, m_1\}, \dots, \min\{n_r, m_r\}), \\
& \max \{\n,\m\} = (\max \{n_1, m_1\}, \dots, \{n_r, m_r\}).
\end{align*}
This leads naturally to the following definition of $\v(M)$ which we
shall often make use of throughout the paper.

\begin{defn}\label{def.v} Suppose $M$ is minimally generated in degrees $\d_1, \dots, \d_u \in \ZZ^r$, then
we define
$$\v(M) = \min \{ \d_1, \dots, \d_u \}.$$
\end{defn}

Besides the given $\NN^r$-graded structure, $S$ has a natural
$\NN$-graded structure defined by $S = \bigoplus_{n \in \ZZ} S_n$
where $S_n = \bigoplus_{|\n| = n} S_\n$ and $|\n|$ indicates the
sum of all components in $\n$. Let $\fM_S$ be the maximal
homogeneous ideal of $S$ with respect to this grading. Observe that
$\fM_S$ is also $\NN^r$-homogeneous. Thus, the local cohomology
modules, $H^i_{\fM_S}(M)$, are $\ZZ^r$-graded modules for all $i$.
This leads to a natural multigraded analog of the usual $\ZZ$-graded
$a$-invariant. The multigraded $a$-invariant is well-defined and has
been studied in more detail in \cite{ha-reg, HHR1, HHR2}.

\begin{defn} \label{def.a invariant}
For each $j = 1, \dots, r$, let $$a_j(M) = \max \big\{ m \in \ZZ
~\big|~ [H^{\dim M}_{\fM_S}(M)]_\n \not= 0 \mbox{ for some } \n \in \ZZ^r \mbox{ with } n_j = m\big\}.$$ The {\it
multigraded a-invariant} of $M$ is define to be the vector $$\mba(M) = (a_1(M),
\dots, a_r(M)) \in \ZZ^r.$$
When $r = 1$, we shall omit the vector notation and simply consider $\mba(M)$ as an integer.
\end{defn}

We shall now recall basic definitions of multi-Rees algebras and modules.

\begin{defn} \label{def.Reesmodule} Let $B$ be a Noetherian ring and let $N$ be a $B$-module.
\begin{enumerate}
\item Let $I \subset B$ be a proper ideal.
The {\it Rees algebra} of $I$ over $B$ is defined to be the subring
$$\R_B(I) = B[It] = \bigoplus_{n \ge 0} I^nt^n \subset B[t].$$ The
{\it Rees module} of $I$ with respect to $N$ is defined to be the
$\R_B(I)$-module $$\R_N(I) = \bigoplus_{n \ge 0}I^nt^nN.$$
\item Let $I_1, \dots, I_r \subset B$ be proper ideals. The {\it multi-Rees algebra} of $I_1, \dots,
I_r$ over $B$ is defined to be the subring
$$\R_B(I_1, \dots, I_r) = \bigoplus_{\n \ge \0}
I_1^{n_1}t_1^{n_1} \dots I_r^{n_r}t_r^{n_r} \subset B[t_1, \dots,
t_r].$$ The {\it multi-Rees module} of $I_1, \dots, I_r$ with
respect to $N$ is defined to be the $R_B(I_1, \dots, I_r)$-module
$$\R_N(I_1, \dots, I_r) = \bigoplus_{\n \ge \0}
I_1^{n_1}t_1^{n_1} \dots I_r^{n_r} t_r^{n_r}N.$$
\end{enumerate}
\end{defn}

Let $Z = \Proj S$ with respect to the $\NN^r$-graded structure of
$S$. As a set, $Z = \{ \pp \in \Spec S ~|~ \pp \mbox{ is }
\NN^r\mbox{-homogeneous and } S_+ \not\subset \pp\}$. The simplest
example is when
%  $S = A[x_{10}, \dots, x_{1N_1}, \dots, x_{r0}, \dots, x_{rN_r}]$
$S = A[x_{ij}]$ where $1 \le i \le r$ and $0 \le j \le N_i$ and $Z =
\PP^{N_1}_A \times \dots \times \PP^{N_r}_A$. Throughout the paper,
let $ R= \R_S(S_+)= S[S_+t]$ be the Rees algebra of $S_+$ over $S$.
Observe that $R$ is a standard $\NN$-graded algebra over $R_0 = S$
where the grading is given by the power of $t$ appearing in each
element. Let $\fM_R$ be the maximal homogeneous ideal of $R$, and
let $Y = \Proj R$ with respect to the $\NN$-graded structure of $R$.
We can view $Y$ as a vector bundle over $Z$ by \cite[Lemma 3.1]{HY}
stated in the next lemma. This provides a natural projection $\pi: Y
\rightarrow Z$.

\begin{lemma} \label{lemma.vectorbundle}
With the above notation, we have
$$Y = \Proj \operatorname{Sym}(\O_Z(\e_1) \oplus \dots \oplus \O_Z(\e_r)).$$
\end{lemma}

%  Needs a better transition from lemma to discussion...

One of the techniques that we employ is to view graded algebras and
modules under various gradings. A simple fact we often use is that
local cohomology modules behave well under a change of grading. More
precisely, suppose $B$ is a standard $\NN^k$-graded algebra over
$(A,\gm)$ and $\aa \subseteq B$ is an $\NN^k$-graded homogeneous
ideal. The local cohomology functors $H^i_\aa(\bullet)$ can be
defined in the category of $\ZZ^k$-graded $B$-modules as usual. That
means if $N$ is a finitely generated $\ZZ^k$-graded $B$-module, then
$H^i_\aa(N)$ is also a $\ZZ^k$-graded $B$-module for all $i$. Let
$\phi: \ZZ^k \rightarrow \ZZ^l$ be a group homomorphism such that
$\phi(\NN^k) \subseteq \NN^l$, and let $B^\phi = \bigoplus_{\m \in
\ZZ^l} \big( \bigoplus_{\phi(\n) = \m} B_\n \big)$ and $N^\phi =
\bigoplus_{\m \in \ZZ^l} \big( \bigoplus_{\phi(\n) = \m} N_\n
\big)$. Then $B^\phi$ is a $\NN^l$-graded ring and $N^\phi$ is a
$\ZZ^l$-graded $B^\phi$-module. It can be seen that
$\big(H^i_\aa(\bullet)\big)^\phi$ and $H^i_{\aa^\phi}(\bullet^\phi)$
are both $\delta$-functors and coincide when $i = 0$. Thus,
$\big(H^i_\aa(N)\big)^\phi = H^i_{\aa^\phi}(N^\phi)$. Hence, when a
new multigraded structure is specified by a group homomorphism
$\phi$, we shall omit the functorial notation $\bullet^\phi$ and
simply write $H^i_\aa(N)$ for $\big(H^i_\aa(N)\big)^\phi$.

A module that we consider under alternate grading is the Rees module
of $S_+$ with respect to M, namely $T = \R_M(S_+)$.
Clearly, $T$ is a $\ZZ$-graded
$R$-module.  Observe that $R$ and $T$ both further possess a
$\ZZ^{r+1}$-graded structure given by
$$R = \bigoplus_{\n \in \NN^r, k \ge 0} R_{(\n;k)} \mbox{ and }
T = \bigoplus_{\n \in \ZZ^r, k \ge 0}T_{(\n;k)},$$ where $R_{(\n;k)}
= S_{(n_1+k, \dots, n_r+k)} t^k \mbox{ and } T_{(\n;k)} =
M_\n[S_+^k]_{(k,\dots,k)}t^k$ for $\n = (n_1, \dots, n_r)$.

The following observation shall prove useful.  Let $W$ be an
arbitrary finitely generated $\ZZ^{r+1}$-graded $R$-module. By
writing $W = \bigoplus_{k \in \ZZ} W_{\bullet; k}$ where
$W_{\bullet; k} = \bigoplus_{\n \in \ZZ^r} W_{\n;k}$, we can
consider $W$ as a $\ZZ$-graded $R$-module. Let $\widetilde{W}$ be
the associated coherent sheaf of $W$ on $Y$. Note the diagonal
subalgebra $S^\Delta \simeq \bigoplus_{k \ge 0} R_{(\0;k)}$, and
there is a canonical isomorphism $Z = \Proj S \simeq \Proj
S^\Delta$. It can be seen that $\pi_* \widetilde{W} = \bigoplus_{\n
\in \ZZ^r} \widetilde{W_{\n; \bullet}}$ where $W_{\n; \bullet} =
\bigoplus_{k \in \ZZ} W_{\n;k}$ is a graded $S^\Delta$-module. The
module $\Gamma(Y, \widetilde{W}) = \bigoplus_{\n \in \ZZ^r}
\Gamma(Z, \widetilde{W_{\n;\bullet}})$ has a natural structure of a
$\ZZ^r$-graded $S$-module. We, therefore, may consider $\Gamma(Y,
\widetilde{\bullet})$ as a functor from the category of
$\ZZ^{r+1}$-graded $R$-modules to the category of $\ZZ^r$-graded
$S$-module.

\begin{lemma} \label{lemma.co-ism}
Let $W$ be a finitely generated $\ZZ^{r+1}$-graded $R$-module, $\M$
be the associated coherent sheaf of $M$ on $Z$, $T = \R_M(S_+)$, and
$\T$ be the associated coherent sheaf of $T$ on $Y$. Then
\begin{enumerate}
\item[(a)] We have isomorphisms
  $$H^i(Y, \widetilde{W}) \simeq H^i(Z, \pi_* \widetilde{W}) \simeq
  \bigoplus_{\n \in \ZZ^r} H^i(Z, \widetilde{W_{\n;\bullet}})
  \ \text{for all} \ i \ge 0.$$
  In particular, for $W = T$, we get
$$H^i(Y, \T) \simeq \bigoplus_{\n \ge \v(M)} H^i(Z, \M(\n)) \ \text{for all} \ i \ge 0.$$
\item[(b)] Let $E = Z \times_A A/\gm$. We have isomorphisms
$$H^i_{\pi^{-1}(E)}(Y,\widetilde{W}) \simeq H^i_E(Z, \pi_* \widetilde{W}) \simeq
\bigoplus_{\n \in \ZZ^r} H^i_E(Z, \widetilde{W_{\n;\bullet}}) \
\text{for all} \ i \ge 0.$$
\end{enumerate}
\end{lemma}

\begin{proof} The first statement of (a) and (b) follow from the arguments
of \cite[p. 322]{HY}. The second statement of (a) takes into account
the canonical isomorphism $Z \simeq \Proj S^\Delta$ and the fact
that $T_{(\n;k)} = M_\n[S_+^k]_{(k,\dots,k)}t^k = 0$ for $\n \not\ge
\v(M)$.
\end{proof}

%%%%%%%%%%%%%%%%%%%%%%%%%%%%%%%%%%%%%%%%%%%%%%%%%%%%%%%%%%%%%%%%%%%%%%%%%%%

\section{Cohen-Macaulay multigraded modules} \label{s.main}

In this section, we prove our first main result. The theorem is stated as follows.

\begin{theorem} \label{mainthm}
Let $S$ be a standard $\NN^r$-graded algebra over a local ring
$(A,\gm)$ such that $S_+$ has positive height, and let $M$ be a
finitely generated $\ZZ^r$-graded $S$-module. Let $Z = \Proj S$ and
let $E = Z \times_A A/\gm$. Let $\M$ be the associated coherent
sheaf of $M$ on $Z$. Then $M$ is a Cohen-Macaulay module with
$\mba(M) < \v(M)$ if and only if the following conditions are
satisfied:
\begin{enumerate}
\item $\Gamma(Z, \M(\n)) = M_\n$ for all $\n \ge \v(M)$,
\item $H^i(Z, \M(\n)) = 0$ for all $i > 0$ and $\n \ge \v(M)$,
\item $H^i_E(Z, \M(\n)) = 0$ for all $i < \dim M - r$ and $\n < \v(M)$.
\end{enumerate}
\end{theorem}

To prove Theorem \ref{mainthm} we shall need some auxiliary results.
As indicated in the introduction, we begin by showing that
$[H^i_{\fM_R}(T)]_{(\n;k)} = 0$ for $i \ge 0, k \ge 0$, and $\n <
\v(M)$ in Lemma \ref{L:3.3} where $T$ denotes the multi-Rees module
$\R_M(S_+)$. Our proof of Lemma \ref{L:3.3} is based upon a simple
observation that if a $\ZZ^l$-graded module $P$ has the
$\ZZ^l$-graded homogeneous decomposition being $P = \bigoplus_{m_1 =
t} P_\m$ for a fixed $m \in \mathbb Z$, then for any $\m \in \ZZ^l$
such that $m_1 \not= t$ we must have $P_\m = 0$.

For any $\ZZ^{r+1}$-graded $R$-module $N = \bigoplus_{\n \in \ZZ^r,
k \in \ZZ} N_{(\n;k)}$, we define the {\it defining region} of $N$
to be
$$\D(N) = \{(\n;k) ~|~ N_{(\n;k)} \not= 0\}.$$

\begin{lemma}\label{L:3.3} For all $i \ge 0$, $k \ge 0$, and
$\n < \v(M)$, we have $[H^i_{\fM_R}(T)]_{(\n;k)} = 0$.
\end{lemma}

\begin{proof} Let $Q$ be the ideal $\bigoplus_{\n > \0} R_{(\n;k)}$ of
  $R$ when $R$ is viewed as a $\NN^{r+1}$-graded ring. Let $R_+$ be
  the irrelevant ideal of $R$ when $R$ is viewed as a $\NN^r$-graded
  ring, i.e. $R_+ = \bigoplus_{(n_1+k, \dots, n_r+k) > \0}
  R_{(\n;k)}$. Then the sequences
$$
0 \rightarrow QT \rightarrow T \rightarrow T/QT \rightarrow 0 \text{
and } 0 \rightarrow R_+T \rightarrow T \rightarrow T/R_+T
\rightarrow 0
%\label{lemma33-eq1}
$$
are exact. By taking the long exact sequences of cohomology, we get
\begin{equation}\label{long exact seq 1}
 \cdots \rightarrow H^{i-1}_{\fM_R}(T/QT) \rightarrow
 H^i_{\fM_R}(QT) \rightarrow H^i_{\fM_R}(T) \rightarrow
 H^i_{\fM_R}(T/QT) \rightarrow \cdots
\end{equation}
and
\begin{equation}\label{long exact seq 2}
 \cdots \rightarrow H^{i-1}_{\fM_R}(T/R_+T) \rightarrow
 H^i_{\fM_R}(R_+T) \rightarrow H^i_{\fM_R}(T) \rightarrow
 H^i_{\fM_R}(T/R_+T) \rightarrow \cdots.
\end{equation}

Let $\d_i=(d_{i,1}, \cdots, d_{i,r})$, $1 \leq
i \leq u$ be the degree of a minimal generator of $M$ as in
Definition~\ref{def.v}.
For a region $D \subseteq \ZZ^{r+1}$, we shall denote
$\bigoplus_{(\n;k) \in D} T_{(\n;k)}$ by $T_D$. Observe that the
defining region of $T$ is $\D(T) = \{ (\n;k) ~|~ k \ge 0 \text{ and
} \exists j: \n \ge \d_j \},$ and the defining region of $QT$ is
$\D(QT) = \{ (\n;k) ~|~ k \ge 0 \text{ and } \exists j: \n > \d_j
\}$. It is also easy to see that $[T/QT]_{(\n;k)} \not= 0$ only if
$(\n;k) \in \D(T) \backslash \D(QT) = \{ (\n;k) ~|~ k \ge 0 \text{
and } \exists j, l: \n \ge \d_j, n_l = d_{j,l} \}.$

Suppose $\v(M) = (v_1, \dots, v_r)$.
We define $d_{\max} = \max \{ d_{i,r} ~|~ 1
\le i \le u\}$ and let $[r-1] = \{1, \dots, r-1\}$. For $1 \le j \le
u$ and a nonempty index set $I \subset [r-1]$, let
$$
A_{j,I} = \{ (\n;k) ~|~ k \ge 0, n_w > d_{j,w} \ \forall w \not\in
I, n_s = d_{j,s} \ \forall s \in I, \mbox{ and } n_r > d_{\max} \}.
$$
For $v_r \le t \le d_{\max}$, let $$B_t = \{ (\n;k) ~|~ k \ge 0, n_r
= t\}.$$ It can be seen that the regions $\{A_{j,I}, B_t ~|~ 1 \le j
\le u, \emptyset \not= I \subset [r-1], v_r \le t \le d_{\max}\}$
are pairwise disjoint.  It also follows from the definition of
$A_{j,I}$'s and $B_t$'s that
$$
\D(T) \backslash \D(QT) \subseteq \big(\bigcup_{j=1}^u
\bigcup_{\emptyset \not= I \subset [r-1]} A_{j,I}\big) \bigcup
\big(\bigcup_{t=v_r}^{d_{\max}} B_t \big).
$$

Thus, we can write
$$
T/QT = \big(\bigoplus_{j=1}^u \bigoplus_{\emptyset \not= I \subset
[r-1]} [T/QT]_{A_{j,I}} \big) \bigoplus
(\bigoplus_{t=v_r}^{d_{\max}} [T/QT]_{B_t}\big).
$$

The importance of the regions $A_{j,I}$'s and $B_t$'s lies in the
fact that $[T/QT]_{A_{j,I}}$ and $[T/QT]_{B_t}$ are submodules of
$T/QT$ for all $j$, $I$, and $t$. Since local cohomology commutes
with direct sum, this allows us to get the following decomposition
of $H^i_{\fM_R}(T/QT)$ into a direct sum of submodules defined over
the $A_{j,I}$'s and $B_t$'s:
\begin{align}
H^i_{\fM_R}(T/QT) = \big(\bigoplus_{j=1}^u \bigoplus_{\emptyset
\not= I \subset [r-1]} H^i_{\fM_R}([T/QT]_{A_{j,I}})\big) \bigoplus
\big( \bigoplus_{t=v_r}^{d_{\max}} H^i_{\fM_R}([T/QT]_{B_t})\big).
\label{lemma33-eq21}
\end{align}

Observe that $[T/QT]_{A_{j,I}}$ is annihilated by $\bigoplus_{n_s >
0 \ \forall s \in I} R_{(\n;k)}$ in $T/QT$, and so
$[T/QT]_{A_{j,I}}$ can be viewed as a $\ZZ^{r+1-|I|}$-graded module
over $R_I = \bigoplus_{n_s=0 \ \forall s \in I} R_{(\n; k)}$. It now
follows from the definition of local cohomology that
$$H^i_{\fM_R}([T/QT]_{A_{j,I}}) = H^i_{\fM_{R_I}}([T/QT]_{A_{j,I}})$$
is a $\ZZ^{r+1-|I|}$-graded module over $R_I$. Moreover,
$H^i_{\fM_{R_I}}([T/QT]_{A_{j,I}})$ has the following
$\ZZ^{r+1-|I|}$-graded decomposition
$$H^i_{\fM_{R_I}}([T/QT]_{A_{j,I}}) = \bigoplus_{n_s = d_{j,s} \ \forall s \in I} H^i_{\fM_{R_I}}([T/QT]_{A_{j,I}})_{(\n;k)}.$$
This implies that if $n_s \not= d_{j,s}$ for some $s \in I$ then the
term $H^i_{\fM_R}([T/QT]_{A_{j,I}})_{(\n;k)}$ is not present in the
homogeneous decomposition of $H^i_{\fM_R}([T/QT]_{A_{j,I}})$. That
is, for $\n \in \ZZ^r$ such that $n_s \not= d_{j,s}$ for some $s \in
I$, we must have $H^i_{\fM_R}([T/QT]_{A_{j,I}})_{(\n;k)} = 0.$

By a similar argument, we have $H^i_{\fM_R}([T/QT]_{B_t})_{(\n;k)} =
0$ if $n_r \not= t$. Hence, it follows from (\ref{lemma33-eq21})
that
$$H^i_{\fM_R}(T/QT)_{(\n;k)} = 0  \text{ if } n_s \not= d_{j,s}  \text{ for all } 1 \le s \le r-1 \mbox{ and } n_r < v_r.$$
This, together with the definition of $\v(M)$, implies that
\begin{align}
H^i_{\fM_R}(T/QT)_{(\n;k)} = 0 \text{ for all } \n < \v(M). \label{lemma33-eq4}
\end{align}

By a similar line of arguments on the defining regions of $T$ and $R_+T$, we have
\begin{align}
H^i_{\fM_R}(T/R_+T)_{(\n,k)} = 0 \text{ for all } k > 0. \label{lemma33-eq5}
\end{align}

Observe further that there is an obvious isomorphism $Q \rightarrow
R_+(-\1,1)$ which maps $R_{(\n;k)}$ to $R_{(\n-\1;k+1)}$. Hence, it
follows from (\ref{long exact seq 1}), (\ref{long exact seq 2}),
(\ref{lemma33-eq4}) and (\ref{lemma33-eq5}) that, for any $k \ge 0$
and $\n < \v(M)$,
\begin{align}
H^i_{\fM_R}(T)_{(\n;k)} \simeq H^i_{\fM_R}(QT)_{(\n;k)} \simeq
H^i_{\fM_R}(R_+T)_{(\n-\1;k+1)} \simeq H^i_{\fM_R}(T)_{(\n-\1;k+1)}.
\label{lemma33-iso}
\end{align}
Moreover, $H^i_{\fM_R}(T)_{(\n;k)} = 0 $ for $k \gg 0$. Therefore,
by successively applying (\ref{lemma33-iso}), we have
$H^i_{\fM_R}(T)_{(\n;k)} = 0$ for all $\n < \v(M)$. The lemma is
proved.
\end{proof}

Lemma \ref{L:3.3} and the sequence (\ref{SSsequence}) imply that
$[H^i_{\fM_S}(M)]_\n = [H^i_F(Y,\T)]_\n$ for $\n < \v(M)$. Thus, to
characterize the Cohen-Macaulayness of $M$, or equivalently, the
vanishing of $[H^i_{\fM_S}(M)]_\n$ for $i < \dim M$ and $\n \in
\ZZ^r$, we proceed by relating $H^i_F(Y, \T)$ to $H^{i-r}_E(Z,\M)$
and establishing the vanishing of $[H^i_{\fM_S}(M)]_\n$ for $\n
\not< \v(M)$.

\begin{lemma} \label{P:3.1}
Let $E=Z \times_A A/\gm$ and $F= Y \times_S S/\fM_S$ where $\fM_S$ is
the homogeneous maximal ideal of $S$.  Then as a $\ZZ^r$-graded $S$-module
$$H^i_F(Y,\T) = \bigoplus_{\n < \v(M)} H^{i-r}_E(Z,\M(\n)) \text{ for all } i \geq 0. $$
\end{lemma}

\begin{proof}
Let $G= Y \times_S S/S_*$ where $S_* = \bigoplus_{\n \neq \0}
S_{\n}$.  As noted in the preliminaries, we consider the functor
$\Gamma_F(Y, \tilde{\bullet})$ from the category of
$\ZZ^{r+1}$-graded $\R_S(S_+)$-modules to the category of
$\ZZ^r$-graded $S$-modules.  Since $\fM_S = \gm \oplus S_*$, this
functor is equal to the composition functor
$\Gamma_{\pi^{-1}(E)}(Y,\H^0_G(\tilde{\bullet}))$.  It follows that
there is a spectral sequence
\begin{align}
E^{p,q}_2 = H^p_{\pi^{-1}(E)}(Y,\H^q_G(\T)) \Rightarrow H^{p+q}_F(Y,\T). \label{specseq}
\end{align}
On the other hand, by Lemma \ref{lemma.co-ism}
$$
H^p_{\pi^{-1}(E)}(Y,\H^q_G(\T)) = H^p_E(Z, \pi_*(\H^q_G(\T)))
$$
as $\ZZ^r$-graded $S$-modules. Now, the conclusion follows from
Lemma \ref{L:3.2} which shows that the spectral sequence
(\ref{specseq}) degenerates.
\end{proof}

\begin{lemma} \label{L:3.2}
Let $G=Y \times_S S/S_*$ where $S_* = \bigoplus_{\n \neq \0}
S_{\n}$.  Then $\H^i_G(\T) = (H^i_{S_*R}(T))^{\sim}$. Moreover, if
$Z=\Proj S$ and $\pi: Y \to Z$ is the canonical projection, we have
$$ \pi_*(\H^i_G(\T)) = \left\{
\begin{array}{ll}
    0, & \text{if } i \neq r, \\
    \bigoplus_{\n < \v(M)} \M(\n), & \text{if } i = r \\
\end{array}
\right. $$
as $\ZZ^r$-graded $\O_Z$-modules.
\end{lemma}

\begin{proof}
For any affine open set $D_+(f) \subset Y$ where $f \in R$ is a homogeneous element, we have
$$
\H^i_G(\T)\big|_{D_+(f)} = (H^i_{S_*R_{(f)}}(T_{(f)}))^{\sim} = ((H^i_{S_*R}(T))_{(f)})^{\sim}.
$$
This proves the first claim.

To prove the second claim, we consider the shifted module $N =
M(\v(M))$ and let $W = \R_N(S_+)$ be the Rees module of $S_+$ with
respect to $N$. As $T$, $W$ admits a $\ZZ^{r+1}$-graded structure $W =
\bigoplus_{\n \in \ZZ^r, k \ge 0} W_{(\n;k)}$ and there is a natural
isomorphism $W \rightarrow T(\v(M), 0)$. It then follows from Lemma
\ref{lemma.co-ism} and the preceding discussion that $\pi_*
(\H^i_G(\widetilde{W})) = \pi_* (\H^i_G(\T)) \otimes \O_Z(\v(M))$.

Now cover $Z$ with open affine sets $\{U= \Spec S_{(s_1 \cdots s_r)}
~|~ s_j \in S_{\e_j} \ \forall j = 1, \ldots, r\}$. By construction,
$\pi^{-1}(U) = \Spec R_{(s_1 \cdots s_r t)}$. Notice that $W_k= N
(S_+)^k t^k$. So $W_{(s_1\cdots s_rt)}$ contains elements in the
form of $\frac{mft^k}{(s_1\cdots s_rt)^k}$ with $m \in N$ and $f \in
(S_+)^k = \bigoplus_{\n \geq \0}S_{(n_1+k, \dots, n_r+k)}$. Since $N
= M(\v(M))$, we have $\deg(m) \geq \0$. Therefore, we can write
$\frac{mf}{(s_1\cdots s_r)^k}$ as a sum of forms like
$h(\frac{s_1}{1})^{\ell_1} \cdots (\frac{s_r}{1})^{\ell_r}$ with $h
\in N_{(s_1\cdots s_r)}$ and $\ell_i \geq 0$. Set $B = N_{(s_1
\cdots s_r)}$ and $t_j = s_j/1 \in W_{(s_1 \cdots s_r t)}$ for $j =
1, \ldots, r$. From the above observation, $W_{(s_1 \cdots s_r t)} =
B[t_1, \dots, t_r]$ as $S_{(s_1 \cdots s_r)}$-module. Since $G \cap
\pi^{-1}(U) = \V(t_1, \dots, t_r)$, we have
$$\H^i_G(\widetilde{W})\big|_{\pi^{-1}(U)} = H^i_{(t_1, \dots, t_r)}(B[t_1, \dots, t_r])^{\sim}.$$
Moreover, it follows from \cite[Remarque 2.1.11 of Chapitre
III]{EGA3} that
$$H^i_{(t_1, \dots, t_r)}(B[t_1, \dots, t_r]) = \left\{
\begin{array}{ll}
    0, & \text{if } i \neq r, \\
    \bigoplus_{\n < \0} B t_1^{n_1} \cdots t_r^{n_r}, & \text{if } i = r. \\
\end{array}
\right.$$
Thus,
$$\pi_* (\H^i_G(\widetilde{W})) = \left\{
\begin{array}{ll}
    0, & \text{if } i \neq r, \\
    \bigoplus_{\n < \0} \widetilde{N}(\n), & \text{if } i = r. \\
\end{array}
\right. $$
The proof is completed by observing that $\widetilde{N} = \M(\v(M)).$
\end{proof}

We are now ready to prove Theorem \ref{mainthm}.

\noindent\textbf{Proof of Theorem \ref{mainthm}.} We first observe
that the Sancho de Salas sequence (\ref{SSsequence}) is an exact
sequence of $\ZZ^r$-graded modules.

It can be seen that $M$ is a Cohen-Macaulay module with $\mba(M) <
\v(M)$ if and only if the following conditions are satisfied:
\begin{enumerate}
\item[(i)] $H^i_{\fM_S}(M)_\n = 0$ for all $0 \le i < \dim M$ and $\n < \v(M)$,
\item[(ii)] $H^i_{\fM_S}(M)_\n = 0$ for all $i \ge 0$ and $\n \not< \v(M)$.
\end{enumerate}

It follows from Lemma \ref{L:3.3} that $[H^i_{\fM_R}(T)]_{(\n;0)} =0$
for all $\n < \v(M)$.  By the exact sequence (\ref{SSsequence}), this
shows that $[H^i_{\fM_S}(M)]_{\n} = [H^i_F(Y, \T)]_{\n}$ for all $\n <
\v(M)$ and all $i \ge 0$.  Using Lemma~\ref{P:3.1}, it follows that
\begin{align}
& H^i_{\fM_S}(M)_{\n} = H^{i-r}_E(Z,\M(\n)) \text{ for all } i \ge 0
\text{ and } \n < \v(M), \label{main theorem eq1}
\\ & H^i_{\fM_S}(M)_\n = 0 \text{ for all } i \geq 0 \text{ and }
  \n \not< \v(M). \label{main theorem eq2}
\end{align}
Thus, $[H^i_{\fM_S}(M)]_{\n} = 0$ for $i < \dim M$ and $\n < \v(M)$
if and only if (3) holds. That is, (i) is equivalent to (3).

On the other hand using Lemma \ref{P:3.1} and the Sancho de Salas
sequence (\ref{SSsequence}), we have that $[H^i_{\fM_R}(T)]_0 \simeq
\bigoplus_{\n \not< \v(M)} [H^i_{\fM_S}(M)]_\n$.
%This and (\ref{main theorem eq2}) implies that
%$$[H^i_{\fM_R}(T)]_0 = \bigoplus_{\n \not<
%\v(M)}[H^i_{\fM_S}(M)]_{\n}.$$
By (\ref{main theorem eq2}), (ii) is equivalent to the condition
that $[H^i_{\fM_R}(T)]_{0}=0$. It then follows from the
Serre-Grothendieck correspondence between local cohomology and sheaf
cohomology that this is equivalent to having $\Gamma(Y,\T) \simeq
T_0 = M$ and $H^i(Y,\T) = 0$ for $i >0$. Moreover, by Lemma
\ref{lemma.co-ism}, $\Gamma(Y,\T) \simeq \bigoplus_{\n \ge \v(M)}
\Gamma(Y, \M(\n))$ and $H^i(Y,\M) \simeq \bigoplus_{\n \ge
\v(M)}H^i(Y,\M(\n))$. Thus, (ii) is equivalent to (1) and (2). \qed

%%%%%%%%%%%%%%%%%%%%%%%%%%%%%%%%%%%%%%%%%%%%%%%%%%%%%%%%%%%%%%%%%%%%%%%

\section{Cohen-Macaulay multi-Rees modules} \label{s.applications}

In this section we shall apply our main result, Theorem
\ref{mainthm}, to investigate the Cohen-Macaulayness of multi-Rees
modules. Our work extends previous studies on the Cohen-Macaulayness
of multi-Rees algebras in \cite{HY}.

Throughout this section, $(A,\gm)$ is a local ring, $N$ is a
finitely generated $A$-module, $I_1, \dots, I_r \subset A$ are
ideals of positive heights with respect to $N$, and $S = \R_A(I_1,
\dots, I_r)$ is the multi-Rees algebra of $I_1, \dots, I_r$.
Clearly, $S$ is a standard $\NN^r$-graded algebra over $A$. Observe
further that the Rees algebra $\R_A(I_1 \cdots I_r)$ is the diagonal
subalgebra $S^\Delta$ of $S$. Let $M = \R_N(I_1, \dots, I_r)$ be the
multi-Rees module of $I_1, \dots, I_r$ with respect to $N$ as
defined in Section~\ref{s.prel}. Then $M$ is a finitely generated
$\ZZ^r$-graded $S$-module, and similarly, the Rees module $\R_N(I_1
\cdots I_r)$ is the diagonal submodule $M^\Delta$ of $M$.

\begin{lemma}\label{L:2.1} With notation as above, we have $\mba(M) = -\1$.
\end{lemma}

\begin{proof} We shall use induction on $r$. For $r = 1$, our argument
is similar to that of \cite[Lemma 2.1]{HHR1}. Observe first that
when $r = 1$, $S = A[I_1t]$ is a standard $\NN$-graded over $A$ and
$M = \R_N(I_1)$ is a finitely generated $\ZZ$-graded $S$-module. Let
$S_+$ be the homogeneous irrelevant ideal of $S$ under this grading
(i.e. $S_+ = (I_1t)S$). It can be seen that $\fM_S = \gm \oplus S_+$
where $\fM_S$ is the maximal ideal of $S$. Let $\G_N = M/I_1M \simeq
\bigoplus_{k \ge 0} I_1^kN/I_1^{k+1}N$. We have the following exact
sequences
\begin{equation*}
0 \rightarrow S_+M \rightarrow M \rightarrow M/S_+M \simeq N
\rightarrow 0
\end{equation*}
\begin{equation*}
0 \rightarrow S_+M(1) \rightarrow M \rightarrow
\G_N \rightarrow 0.
\end{equation*}
taking the corresponding long exact sequences of local cohomology
modules, for any $i \ge 0$ and $n \in \ZZ$, we have
\begin{align}
& \cdots \rightarrow [H^{i-1}_{\fM_S}(N)]_n \rightarrow
[H^i_{\fM_S}(S_+M)]_n \rightarrow [H^i_{\fM_S}(M)]_n \rightarrow
[H^i_{\fM_S}(N)]_n \rightarrow \cdots \label{lemma1.eq1} \end{align}
and
\begin{align}
& \cdots \rightarrow [H^{i-1}_{\fM_S}(\G_N)]_n \rightarrow
[H^i_{\fM_S}(S_+M)]_{n+1} \rightarrow [H^i_{\fM_S}(M)]_n \rightarrow
[H^i_{\fM_S}(\G_N)]_n \rightarrow \cdots. \label{lemma1.eq2}
\end{align}

Observe that $N \simeq M/S_+M$, as an $S$-module, is concentrated in
degree 0. Thus, for $n \not= 0$, $[H^i_{\fM_S}(N)]_n$ vanish for all
$i \ge 0$. The long exact sequence (\ref{lemma1.eq1}) implies that
\begin{align}
[H^i_{\fM_S}(S_+M)]_n \simeq [H^i_{\fM_S}(M)]_n \mbox{ for any } i \ge
0 \mbox{ and } n \not= 0. \label{lemma1.eq3}
\end{align}
By \cite[Theorem 4.4.6]{BH}, $\dim \G_N = \dim N = \dim \R_N(I_1) -
1 = \dim M - 1$. Thus $H^{\dim M}_{\fM_S}(\G_N) = 0$. Therefore
(\ref{lemma1.eq2}) implies that for any $n \in \ZZ$, there is a
surjection
\begin{align}
[H^{\dim M}_{\fM_S}(S_+M)]_n \smap [H^{\dim M}_{\fM_S}(M)]_n. \label{lemma1.eq4}
\end{align}
It follows from (\ref{lemma1.eq3}) and (\ref{lemma1.eq4}) that for
any $n \not= 0$, there is a surjection
\begin{align}
[H^{\dim M}_{\fM_S}(M)]_n \smap [H^{\dim M}_{\fM_S}(M)]_{n-1}. \label{lemma1.eq5}
\end{align}

By successively applying (\ref{lemma1.eq5}) and noting that $H^{\dim
M}_{\fM_S}(M)_n = 0$ for $n \gg 0$, we have $[H^{\dim
M}_{\fM_S}(M)]_n = 0$ for all $n \ge 0$. We must also have $[H^{\dim
M}_{\fM_S}(M)]_{-1} \not= 0$; otherwise, again by successively
applying (\ref{lemma1.eq5}) we would conclude that $H^{\dim
M}_{\fM_S}(M) = 0$, which is impossible. Hence, $\mba(M) = -1$. The
statement holds for $r = 1$.

Assume that the statement already holds for $r-1$ (for some $r \ge
2$). Let $A' = \R_A(I_r)$. Then $S = A[I_1t_1, \dots, I_rt_r] =
A'[I_1t_1, \dots, I_{r-1}t_{r-1}] = \R_{A'}(I_1, \dots, I_{r-1})$
can be viewed as an $\NN^{r-1}$-graded algebra over $A'$, where the
grading is given by powers of $t_1, \dots, t_{r-1}$. By giving $t_r$
degree 0, we can also view $M$ as a $\ZZ^{r-1}$-graded $S$-module.
Since local cohomology modules behave well under a change of
grading, we have the following $\ZZ^{r-1}$-graded homogeneous
decomposition of $H^i_{\fM_S}(M)$:
$$H^i_{\fM_S}(M) = \bigoplus_{\n' = (n_1, \dots, n_{r-1}) \in \ZZ^{r-1}} [H^i_{\fM_S}(M)]_{\n'},$$
where $[H^i_{\fM_S}(M)]_{\n'} = \bigoplus_{n \in \ZZ}
[H^i_{\fM_S}(M)]_{(n_1, \dots, n_{r-1}, n)}$. By induction, as a
$\ZZ^{r-1}$-graded module, $\mba(M) = -\1 \in \ZZ^{r-1}$. Thus,
$a_j(M) = -1$ for all $j =1, \dots, r-1$. It remains to show that
$a_r(M) = -1$.

Let $A'' = \R_A(I_1, \dots, I_{r-1})$.  Then $S = \R_{A''}(I_r)$ can
now be viewed as an $\NN$-graded algebra over $A''$. By a similar
argument to last paragraph, we can view $M$ as a $\ZZ$-graded
$S$-module; and therefore, by induction, $a_r(M) = -1$. Hence,
$\mba(M) = (a_1(M), \dots, a_r(M)) = -\1 \in \ZZ^r$.
\end{proof}

The next theorem extends a well-known result for multi-Rees algebras to arbitrary multi-Rees modules.

\begin{theorem} \label{Rees module}
Let $(A,\gm)$ be a local ring and let $N$ be a finitely generated
$A$-module. Let $I_1, \dots, I_r \subset A$ be ideals of positive
heights with respect to $N$. Assume that the multi-Rees module
$\R_N(I_1, \dots, I_r)$ is Cohen-Macaulay. Then the Rees module
$\R_N(I_1 \cdots I_r)$ is also Cohen-Macaulay.
\end{theorem}

\begin{proof} Recall that the multi-Rees algebra
$S = \R_A(I_1, \dots, I_r)$ is a standard $\NN^r$-graded algebra
over $A$, and the Rees algebra $\R_A(I_1 \cdots I_r)$ is the
diagonal subalgebra $S^\Delta$ of $S$. Let $Z = \Proj \R_A(I_1,
\dots, I_r)$, and let $E = Z \times_A A/\gm$. As before, there is a
canonical isomorphism $f: Z \stackrel{\simeq}{\longrightarrow} \Proj
S^\Delta = \Proj \R_A(I_1 \cdots I_r)$ given by the inclusion
$S^\Delta \hookrightarrow S$.

Again, let $M = \R_N(I_1, \dots, I_r)$. For simplicity, we denote
$ M^\Delta = \R_N(I_1 \cdots I_r)$ by $L$. By \cite[Theorem 4.4.6]{BH},
we have
\begin{align}
\dim M = \dim N + r = \dim L + (r-1). \label{lemma21.eq1}
\end{align}
Let $\M$ and $\L$ be the associated coherent sheaves of $M$ and $L$
on
 $Z$ and $\Proj S^\Delta$ respectively. It can be seen that $f_*
\M = \L$ and $f^*(\O_{\Proj S^\Delta}(n)) = \O_Z(n, \dots, n)$.
Thus, by the projection formula we get
\begin{align}
f_* \M(n, \dots, n) = f_* (\M \otimes \O_Z(n, \dots, n)) = \L \otimes \O_{\Proj S^\Delta}(n) = \L(n). \label{lemma21.eq0}
\end{align}

It follows from Lemma \ref{L:2.1} that $\mba(M) = -\1 < \0 = \v(M).$
Theorem \ref{mainthm} together with (\ref{lemma21.eq1}) and
(\ref{lemma21.eq0}) imply that
\begin{enumerate}
\item[(i)] $\Gamma(\Proj S^\Delta, \L(n)) = M_{(n, \dots, n)} = L_n$ for all $n \ge 0$,
\item[(ii)] $H^i(\Proj S^\Delta, \L(n)) = 0$ for all $i > 0$ and $n \ge 0$,
\item[(iii)] $H^i_E(\Proj S^\Delta, \L(n)) = 0$ for all $i < \dim L-1$ and $n < 0$.
\end{enumerate}
The Rees module $\R_N(I_1 \cdots I_r)$ is Cohen-Macaulay now follows
from a special use of Theorem~\ref{mainthm} when $r=1$.
\end{proof}

The converse of Theorem~\ref{Rees module} is not always true. The
rest of the paper is devoted to show that when $\Proj S^\Delta$ is a
Cohen-Macaulay scheme, $N$ is free in the punctured
spectrum of $A$, and the analytic spread of $I_1 \cdots I_r$ is small,
there are conditions which, together with the Cohen-Macaulayness of
the usual Rees module $\R_N(I_1 \cdots I_r)$, imply the
Cohen-Macaulayness of the multi-Rees module $\R_N(I_1, \dots, I_r)$.
We recall that $N$ is free in the punctured spectrum of $A$
means $N$ is a free module if localized at every prime ideal in $\Spec A$
with the only possible exception at the maximal ideal. This condition
implies that $\R_N(I_1 \cdots I_r)$ is associated to a locally free
sheaf on $\Proj S^\Delta$.

We shall need some preliminary results. The following lemmas are
generalization of \cite[Lemmas 4.4 and 4.5]{HY} from Rees algebras
to Rees modules. We shall sketch the proof of Lemma~\ref{L:4.4} and
leave that of Lemma~\ref{L:4.5} to the reader.

\begin{lemma}\label{L:4.4}
Assume $(A,\gm)$ is a local ring and $N$ is a finitely generated
$A$-module of dimension $d$. Let $I \subset A$ be an ideal of
positive grade with respect to $N$. Let $P = \R_A(I)$ and $L=
\R_N(I)$.  Let $\L$ be the associated coherent sheaf of $L$ over $Z
= \Proj P$. Let $E = Z \times_A A/\gm$, and let $\ell = \ell(I)$ be
the analytic spread of $I$. Assume that $L$ is Cohen-Macaulay. Then
\begin{enumerate}
\item[(a)] $H^i(Z, \L(\ell - 1 - i)) = 0$
 for all $i >0$;
\item[(b)] If $Z$ is a Cohen-Macaulay scheme and $\L$ is a locally free sheaf, then
$$H^i_E(Z,\L(d - \ell - i)) =0 \mbox{ for all } i <d.$$
\end{enumerate}
\end{lemma}

\begin{proof} For $i > 0$, by the Serre-Grothendieck correspondence we have
$$H^i(Z,\L(\ell-1-i)) \simeq [H^{i+1}_{P_+}(L)]_{\ell-1-i}.$$
By Lemma \ref{L:2.1}, $a(L) = -1$. This fact and the assumption that
$L$ is Cohen-Macaulay imply that $[H^{i+1}_{P_+}(L)]_n = 0$ for all
$i \ge 0$ and $n \ge 0$. Hence, $H^i(Z,\L(\ell-1-i)) = 0$ for all $0
< i \le \ell-1$.

Observe further that by definition, the closed fiber of the
canonical projection $Z \longrightarrow \Spec A$ has dimension
$\ell-1$. Thus, it follows (cf. \cite[Corollaire 4.2.2 in Chapitre
III]{EGA3}) that $H^i(Z,\F) = 0$ for every coherent sheaf $\F$ on
$Z$ if $i \ge \ell$. In particular, this implies that
$H^i(Z,\L(\ell-1-i)) = 0$ for all $i \ge \ell$. We have proved (a).

To prove (b), we first observe that $\v(L) = 0$. By Theorem
\ref{mainthm},
$$H^i_E(Z, \L(d-\ell-i)) = 0 \mbox{ if } i > d-\ell.$$

On the other hand, by Lipman's global-local duality theorem (cf.
\cite[Theorem on p. 188]{LI2}), $H^i_E(Z, \L) \simeq
\Hom_A(\Ext^{d-i}(\L, \omega_Z), E_A(A/\gm))$ where $\omega_Z$ is
the dualizing sheaf on $Z$ and $E_A(A/\gm)$ is the injective hull of
$A/\gm$. Since $\L$ is locally free, by Serre duality theorem, we
also have $\Ext^{d-i}(\L,\omega_Z) \simeq H^{d-i}(Z, \L^\vee \otimes
\omega_Z)$ where $\L^\vee = \Hom_{\O_Z}(\L, \O_Z)$. Thus
\[ H^i_E(Z,\L) \simeq \Hom_A(H^{d-i}(Z, \L^\vee \otimes \omega_Z), E_A(A/\gm)). \]
This implies that if $i \le d-\ell$ (i.e. $d-i \ge \ell$) then,
since the closed fiber of the projection $Z \longrightarrow \Spec A$
has dimension $\ell - 1$, we have $H^{d-i}(Z, \L^\vee(-d+\ell+i)
\otimes \omega_Z) = 0$. That is, if $i \le d-\ell$ then
$H^i_E(Z,\L(d-\ell-i)) = 0$. Hence, (b) is proved.
\end{proof}

\begin{lemma}\label{L:4.5}
Let $(A,\gm)$ be a local ring, and let $I_1, \ldots, I_r \subset A$
be ideals of positive grade with respect to a finitely generated
$A$-module $N$. Let $S = \R_A(I_1, \ldots, I_r)$, $Z = \Proj S$, and
$M = \R_N(I_1, \dots, I_r)$. Let $\M$ be the coherent sheaf
associated to $M$ over $Z$. Then
$$\Gamma(Z, \M(\n-\m)) = \Hom_A(I^{m_1}_1 \cdots I^{m_r}_r, \Gamma(Z, \M(\n))) \text{ for all } \n, \m \ge \0.$$
Moreover,
$$\Gamma(Z, \M(\n-\m)) = \Gamma(Z, \M(\n)) :_{\Gamma(Z, \M)} (I^{m_1}_1 \cdots I^{m_r}_r) \text{ for all } \n \ge \m \ge \0.$$
\end{lemma}

\begin{proof} The proof goes in the same line of arguments as that of \cite[Lemma 4.5]{HY}.
\end{proof}

The next theorem generalizes \cite[Theorem 4.1]{HY} to give the converse of Theorem \ref{Rees module} in the case that $I_1 \cdots I_r$ has small analytic spread.

\begin{theorem} \label{Rees converse}
Let $(A,\gm)$ be a local ring, and let $I_1, \ldots, I_r \subset A$ be
ideals of positive grades with respect to a finitely generated
$A$-module $N$. Assume that $\ell = \ell(I_1 \cdots I_r) \leq 2$.
\begin{enumerate}
\item[(a)] If $\R_N(I_1, \ldots, I_r)$ is Cohen-Macaulay,
then $\R_N(I_1 \cdots I_r)$ is Cohen-Macaulay and the condition
$(I_{j_1} \cdots I_{j_k})N :_N I_{j_l} = (I_{j_1} \cdots I_{j_{l-1}}
  \cdot I_{j_{l+1}} \cdots I_{j_k})N$ holds for all
  $1 \le j_1 < \cdots < j_k \le r$ and $1 \leq l \le k$.
\item[(b)] Conversely, if $\R_N(I_1 \cdots I_r)$ is Cohen-Macaulay, the condition
$(I_{j_1} \cdots I_{j_k})N :_N I_{j_l} = (I_{j_1} \cdots I_{j_{l-1}}
\cdot I_{j_{l+1}} \cdots I_{j_k})N$ holds for all $1 \le j_1 < \cdots
< j_k \le r$ and $1 \leq l \le k$, and if, in addition, $N$ is free
in the punctured spectrum and $\Proj \R_A(I_1 \cdots I_r)$ is a Cohen-Macaulay scheme, then $\R_N(I_1, \dots, I_r)$ is Cohen-Macaulay.
\end{enumerate}
\end{theorem}

\begin{proof}
As before, let $S = \R_A(I_1, \dots, I_r)$, $M = \R_N(I_1, \dots,
I_n)$, $Z = \Proj S$ and $E = Z \times_A A/\gm$.  Let $\M$ be
the associated coherent sheaf of $M$ on $Z$. The first part of (a)
follows from Theorem \ref{Rees module}. To prove the second part of
the statement, we observe that by Theorem \ref{mainthm}, $\Gamma(Z,
\M(\n)) = M_\n = I_1^{n_1} \cdots I_k^{n_k} N$ for all $\n \ge \v(M)
= \0$. Thus, the condition $(I_{j_1} \cdots I_{j_k})N :_N I_{j_l} =
(I_{j_1} \cdots I_{j_{l-1}} \cdot I_{j_{l+1}} \cdots I_{j_k})N$
follows by substituting appropriate $\n$ and $\m$ to Lemma
\ref{L:4.5}.

To prove (b), let $L = M^\Delta = \R_N(I_1 \cdots I_r) $,
and notice that $S^\Delta = \R_A(I_1 \cdots I_r)$ and $\Proj S^\Delta$ is a Cohen-Macaulay scheme. Let $\L$ be the associated coherent sheaf of $L$ on $\Proj S^\Delta$. Recall that
there is a canonical isomorphism $f: Z \rightarrow \Proj S^\Delta$.
Since $L$ is Cohen-Macaulay, by Lemma \ref{L:4.4}(a) and
(\ref{lemma21.eq0}), we have $H^i(Z,\M(\ell - 1 -i,\dots,\ell-1 -i))
= H^i(\Proj S^{\Delta}, \L(\ell-1-i)) = 0$ for $i > 0$. In such a
case, \cite[Lemma 4.2.(a)]{HY} states that $H^i(Z,\M(m_1-i,\dots,m_r
-i)) =0$ if $m_j \geq \ell -1$ for all $j$. Moreover, for $i > 0$,
since $\ell \le 2$, we have $\ell - 1 - i \le 0$. This implies that
\begin{align}
H^i(Z, \M(\n)) = 0 \mbox{ for any } i > 0 \mbox{ and } \n \ge \0. \label{thm41.eq1}
\end{align}
Since $N$ is free in the punctured spectrum of $A$, $\L$ is a locally free sheaf. Thus, similarly, by using Lemma \ref{L:4.4}(b) and \cite[Lemma 4.2(b)]{HY}, we get
\begin{align}
H^i_E(Z, \M(\n)) = 0 \mbox{ for any } i < \dim N \mbox{ and } \n < \0. \label{thm41.eq2}
\end{align}

We shall now verify that
\begin{align}
\Gamma(Z, \M(\n)) = I_1^{n_1} \cdots I_r^{n_r}N \mbox{ for all } \n \ge \0.
 \label{thm41.eq3}
\end{align}
By applying the special case of Theorem \ref{mainthm} for $r = 1$ and
(\ref{lemma21.eq0}), we have
$$\Gamma(Z, \M(m,\dots,m)) = \Gamma(\Proj S^\Delta, \L(m)) = L_m = I_1^m
\cdots I_r^m N \text{ for all } m \ge 0.$$ Now for any $\n \ge \0$,
we can find some $m$ such that $\n \le (m,\dots,m)$. By descending
induction on each coordinate and successively applying Lemma
\ref{L:4.5}, it can be seen that (\ref{thm41.eq3}) holds.

The conclusion of (b) now follows from Theorem \ref{mainthm}
together with (\ref{thm41.eq1}), (\ref{thm41.eq2}), and
(\ref{thm41.eq3}).
\end{proof}

%%%%%%%%%%%%%%%%%%%%%%%%%%%%%%%%%%%%%%%%%%%%%%%%%%%%%%%%%%%%%%%%%%%%%%%%%%%%%%%%

\end{document}